\documentclass[12pt,a4paper]{amsart}
\usepackage{hyperref}
\usepackage{amssymb}
\usepackage{xypic}

\newcommand{\A}{\ensuremath{\mathcal{A}}}

\newcommand{\C}{\ensuremath{\mathcal{C}}}
\newcommand{\D}{\ensuremath{\mathcal{D}}}

\newcommand{\G}{\ensuremath{\mathcal{G}}}

\newcommand{\LLL}{\ensuremath{\mathcal{L}}}
\newcommand{\M}{\ensuremath{\mathcal{M}}}

\newcommand{\PPP}{\ensuremath{\mathcal{P}}}

\newcommand{\T}{\ensuremath{\mathcal{T}}}

\newcommand{\CC}{\ensuremath{\mathbb{C}}}

\newcommand{\NN}{\ensuremath{\mathbb{N}}}
\newcommand{\PP}{\ensuremath{\mathbb{P}}}

\newcommand{\QQ}{\ensuremath{\mathbb{Q}}}

\newcommand{\ZZ}{\ensuremath{\mathbb{Z}}}

\def\hol{{\mathcal{O}}}

\DeclareMathOperator{\Aut}{Aut}
\DeclareMathOperator{\Crit}{Crit}
\DeclareMathOperator{\Ext}{Ext}
\DeclareMathOperator{\Hom}{Hom}

\DeclareMathOperator{\Proj}{\mathbf{Proj}}
\DeclareMathOperator{\Sim}{Sym}
\DeclareMathOperator{\Sing}{Sing}
\DeclareMathOperator{\Span}{Span}
\DeclareMathOperator{\coker}{coker}
\DeclareMathOperator{\rank}{rank}
\DeclareMathOperator{\supp}{supp}

\newtheorem{teo}{Theorem}[section]
\newtheorem{lem}[teo]{Lemma}

\newtheorem{cor}[teo]{Corollary}
\newtheorem{prop}[teo]{Proposition}

\theoremstyle{definition}
\newtheorem{df}[teo]{Definition}

\theoremstyle{remark}
\newtheorem{oss}[teo]{Remark}

\title{Some (big) irreducible components of the moduli space of minimal
   surfaces of general type with $p_g=q=1$ and $K^2=4$}\thanks{Part
   of the article was developped when the author was visiting professor
   at the university of Bayreuth financed by the  {\it DFG
     Forschergruppe ``Klassifikation algebraischer Fl\"achen und
     kompakter komplexer Mannigfaltigkeiten''}\\
2000 Mathematics Subject Classification. Primary 14J29; Secondary 
14D06, 14J10.}
\author{Roberto Pignatelli}
\address{Dipartimento di Matematica, Universit\`a di Trento. Via Sommarive 14,
loc. Povo, I-38050 Trento (Italy)}
\begin{document}

\begin{abstract}
This paper is devoted to the irregular surfaces of general type having
the smallest invariants, $ p_g = q= 1$. We consider the still 
unexplored case $K^2=4$, classifying
those whose Albanese morphism has general fibre of
genus $2$ and such that the direct image
of the bicanonical sheaf under the Albanese morphism is a direct sum 
of line bundles.

We find $8$ unirational families, and we prove that all are irreducible
components of the moduli space of minimal surfaces of general type. This
is unexpected because the assumption on the direct image bicanonical sheaf is
a priori only a closed condition. One more unexpected property is
that all these components have dimension strictly bigger than the expected one.
\end{abstract}

%\thanks{The research of the authors was performed in the realm of the
%DFG SCHWERPUNKT "Globale Methode in der komplexen Geometrie", and of the
%EAGER EEC Project. The third author was  supported by the Schwerpunkt
%and by P.R.I.N. 2002 "Geometria delle variet\`a algebriche" of
%M.I.U.R. and is a member of G.N.S.A.G.A. of I.N.d.A.M. }

\maketitle
\tableofcontents
\pagestyle{myheadings}
%\markboth{$K_S^2=4$, $p_g=q=1$}{$K_S^2=4$, $p_g=q=1$}
\markright{$K_S^2=4$, $p_g=q=1$}

\section*{Introduction}

Minimal surfaces of general type with $p_g=q$ ({\it i.e} with $\chi(\hol)=1$,
the minimal possible value) have attracted the interest of many authors, but we
are very far from a complete classification of them. Bombieri's 
theorem on pluricanonical maps ensures
that there is only a finite number of families of such surfaces, but 
recent results show that
the number of these families is huge (see for instance
  \cite{pragacz},  \cite{bcg}, \cite{bcgp} for the case $p_g=q=0$, 
\cite{PolThesis}, \cite{poliziotto}
for the case $p_g=q=1$, \cite{zucconi} and \cite{penegini} for the 
case $p_g=q=2$).

The irregular case is possibly more affordable, and in fact there is a complete
classification of the case $p_g=q\geq 3$ (\cite{chrisrita}, \cite{pietro}, see
also \cite{survey} for more on surfaces with $\chi(\hol)=1$).

We are interested in the case $p_g=q=1$. A classification of the minimal
surfaces of general type with $p_g=q=1$ and $K^2 \leq 3$ has been obtained
(\cite{bombcat}, \cite{cc1}, \cite{cc2}, \cite{cp}) by looking at the Albanese
morphism,which, for a surface with $q=1$, is a fibration onto an 
elliptic curve.

In this paper we begin the analysis of the next case $K^2=4$, by studying the
surfaces whose general Albanese fibre has the minimal possible genus, 
i.e., genus $2$.

\bigskip

We proved the following
\begin{teo}\label{main}
Let $\M$ be the algebraic subset of the moduli space of minimal 
surfaces of general type
given by the set of isomorphism classes of minimal surfaces $S$ with 
$p_g=q=1$, $K_S^2=4$,
  whose Albanese fibration
$\alpha$ is such that
\begin{itemize}
\item the general fibre of $\alpha$ is a genus $2$ curve;
\item  $\alpha_* \omega_S^2$ is a direct sum of line bundles.
\end{itemize}
Then
\begin{itemize}
\item $\M$ has $8$ connected components, all unirational, one of dimension $5$
and the others of dimension $4$;
\item these are also irreducible components of the moduli space of minimal
surfaces of general type;
\item the general surface in each of these components has ample canonical
class.
\end{itemize}
\end{teo}

We find noteworthy that all these families have bigger dimension than
expected. Standard deformation theory says that any irreducible component of
the moduli space of minimal surfaces of general type containing a surface $S$
has dimension at least $-\chi(\T_S)=10\chi(\hol_S)-2K^2_S$, but by the general
principle ``Hodge theory kills the obstruction'' (stated in \cite{Ran} and
later made precise in \cite{clemens}) this bound is not sharp for irregular
surfaces. By applying this principle as in \cite{cs}, (proof of theorem 5.10),
if $q=1$ a better lower bound is $10\chi(\hol_S)-2K^2_S+p_g=11p_g-2K^2$. This
new bound is sharp, and in fact (\cite{bombcat}, \cite{cc1}, \cite{cc2},
\cite{cp}) all irreducible components of the moduli space of surfaces with
$p_g=q=1$ and $K^2 \leq 3$ attain it. For $K^2=4$ this bound is $3$, and all
our families have strictly bigger dimension.

For technical reasons we assume $\alpha_*\omega_S^2$ to be a sum of line
bundles. This is a closed assumption, and it is rather surprising that all the
families we find are irreducible components of the moduli space of minimal
surfaces of general type. Since \cite{cc1} (thm. 1.4 and prop. 2.2) shows that
the number of direct summands of $\alpha_* \omega_S$ is a topological
invariant, we ask the following

{\bf Question:} is the number of direct summands of $\alpha_* \omega_S^2$ a
deformation or a topological invariant?

The author knows of constructions of minimal surfaces with $p_g=q=1$ 
and $K^2=4$
by Catanese (\cite{cat99}), Polizzi (\cite{poliziotto}) and Rito (\cite{Rito},
\cite{Ritothesis}). Only one of these constructions gives a family of dimension
at least $4$, one of Polizzi's families. But these are obtained by
resolving the singularities of a surface with $4$ nodes; since each of
our $8$ families contains a surface with ample canonical class, the
general surface in each of them is new. In section \ref{def} we show
that the $4$-dimensional family constructed by Polizzi is a proper
subfamily of our ``bigger'' family, the one of dimension $5$.

\bigskip
The proof uses three main tools.

The first one is the study of the relative canonical algebra of genus $2$
fibrations (introduced in \cite{rei} after the results obtained in
  \cite{horP}, \cite{xiao}, \cite{cc1}) and in
particular the structure theorem for genus $2$ fibrations of
\cite{cp}. The assumption on the direct image of the bicanonical 
sheaf is a natural assumption
in view of the results of \cite{cp}.

The second step consists in the analysis of several cases a priori possible:
  some of these are excluded through the investigation of
the geometry of a certain conic bundle, which is obtained  as the 
quotient of our surface by the
involution which induces the hyperelliptic involution on each fibre.
Contradictions are derived  by comparing a ``very negative'' section
$s$ with the branch locus (for example, by showing that $s$ is
contained in its divisorial part, which is reduced, with multiplicity $2$).

Finally, to show that all our families are irreducible components of
the moduli space of minimal surfaces of general type we need to bound
from above  the dimension of the first cohomology group of the tangent
sheaf. To do that, we relate it with the dimension of a subsystem of
the bicanonical system which we can explicitly compute.

\bigskip

The paper is organized as follows.

In section \ref{struct} we recall the structure theorem for genus $2$
fibrations.

In section \ref{8fam}, we apply it to construct $8$ families of minimal
surfaces of general type with $p_g=q=1$, $K^2=4$ whose Albanese fibration
$\alpha$ has fiber of genus $2$ and $\alpha_*\omega_S^2$ is a sum of line
bundles. We also remark that each family contains surfaces with ample canonical
class.

In sections \ref{V1dec} and \ref{V1undec} we show that we have constructed all
  surfaces with the above properties. In other words the image of our families
  in the moduli space of surface of general type equals the scheme $\M$ in
  theorem \ref{main}.

In section \ref{def} we first remark that $\M$ has $8$ unirational connected
components (one for each family) and compute the dimension of each
component. We prove then that they all are irreducible components of the
moduli space of minimal surfaces of general type by investigating
their bicanonical system as mentioned above.

\medskip

{\bf Acknowledgements.} I'm indebted with F.~Catanese for explaining me how to
use the bicanonical curves on a fibred surface to compute its infinitesimal
deformations. I thank F.~Polizzi for many interesting conversations on the
properties of the surfaces in \cite{poliziotto}.

%\section*{Notation}

\section{The structure theorem for genus 2 fibrations}\label{struct}

\subsection{The relative bicanonical map}
In this section we recall results of \cite{cp} (section 4) without
giving any proof. The goal is to explain the structure theorem for
genus $2$ fibration (4.13 there).

Let $f\colon S \rightarrow B$ be a relatively minimal fibration of a smooth
compact complex surface to a smooth curve whose general fibre has genus $2$. We
denote by $F_p$ the fibre $f^{-1}(p)$.

Consider the relative dualizing sheaf $\omega_{S|B}:=\omega_S \otimes
f^{*}\omega_B^{-1}$. The direct images $V_n:=f_* \omega_{S|B}^n$ are vector
bundles on $B$ whose fibre over any point $p$ is canonically isomorphic to
$H^0(\omega_{F_p}^n)$. Therefore the induced rational maps $\varphi_n \colon S
\dashrightarrow \PP(V_n):=\Proj(\Sim V_n)$ (cf. \cite{har}, chapter 2, section
7) map each fibre $F_p$ to the corresponding fibre of $\PP(V_n)$ by its own
n-canonical map.

We remember to the reader that the canonical map of a smooth genus
$2$ curve $F$ is a double cover of $\PP^1$ and that its bicanonical map is the
composition of this map with the $2$-Veronese embedding of $\PP^1$ onto
a conic in $\PP^2$, defined by the isomorphism $\Sim^2(H^0(\omega_F))
\cong H^0(\omega^2_F)$.
The relative analog is an injective morphism of sheaves $\sigma_2
\colon \Sim^2 V_1
\hookrightarrow V_2$ (surjectivity fails on the stalks of points $p$
such that $F_p$ is not $2$-connected) giving a {\em relative
   $2$-Veronese}  $v\colon \PP(V_1) \dashrightarrow \PP(V_2)$ birational
onto a conic bundle $\C$, the image of $\varphi_2$. In fact
$\varphi_2=v \circ \varphi_1$.

The main point is that $\varphi_2$ is always a morphism. More
precisely, $\varphi_2$ is a quasi-finite morphism of degree $2$
contracting exactly the $(-2)$ curves contained in fibres. In other
words, if we substitute $S$ with its relative canonical model (the
surface obtained contracting that curves), $\varphi_2$ becomes a
finite morphism of degree $2$. Moreover $\C$ can only have
singularities of type $A_n$ or $D_n$, that are Rational Double Points.

The structure theorem proves that to reconstruct the pair $(S,f)$ one only
needs to know $\sigma_2$ (that gives at once $\C$ and the isolated
branch points of $\varphi_2$) and the divisorial part $\Delta$ of the
branch locus of $\varphi_2$.
It gives moreover a concrete {\it recipe} to construct all
possible pairs ($\sigma_2,\Delta$).

We now introduce the $5$ {\it ingredients} $(B,V_1,\tau,\xi,w)$,
and then explain how to {\it cook} $\sigma_2$ and $\Delta$ from them.

\subsection{The 5 ingredients}\label{ingr}
Their order is important, since each ingredient is given as an element
in a space that depends on the previously given ingredients. They are
\begin{itemize}
\item[($B$):] Any curve.
\item[($V_1$):] Any rank $2$ vector bundle over $B$.
\item[($\tau$):] Any effective divisor on $B$.
\item[($\xi$):] Any extension class
$$
\xi \in \Ext^1_{\hol_B}(\hol_{\tau},\Sim^2(V_1))/{\Aut_{\hol_B}(\hol_{\tau})}
$$
such that the central term, say $V_2$, of the corresponding short
exact sequence
\begin{equation}\label{conic}
0
\rightarrow
\Sim^2 (V_1)
\stackrel{\sigma_2}{\rightarrow}
V_2
\rightarrow
\hol_{\tau}
\rightarrow
0
\end{equation}
is a vector bundle.
\item[($w$):] A nontrivial element of
$$
\Hom((\det V_1 \otimes \hol_B(\tau))^2,\A_6)/\CC^*.
$$
where $\A_6$ is a vector bundle determined by the other $4$
ingredients as we explain in the following.
\end{itemize}

Consider the map $\nu$ in the natural short exact sequence
$$
0
\rightarrow
(\det V_1)^2
\stackrel{\nu}{\rightarrow}
\Sim^2(\Sim^2(V_1))
\rightarrow
\Sim^4(V_1)
\rightarrow
0;
$$
given locally, if $x_0$, $x_1$ are generators of the stalk of $V_1$ in a
point, by
\begin{equation}\label{nu}
\nu((x_0 \wedge x_1)^{\otimes 2})=(x_0)^2(x_1)^2-(x_0x_1)^2.
\end{equation}
$\A_6$ is the cokernel of the (automatically injective)
composition of maps
\begin{multline}\label{i3def}
(\det V_1)^2 \otimes V_2
\stackrel{\nu \otimes id_{V_2}}\longrightarrow
\Sim^2(\Sim^2(V_1)) \otimes V_2
\stackrel{\Sim^2(\sigma_2)\otimes id_{V_2}}{\longrightarrow}\\
\stackrel{\Sim^2(\sigma_2)\otimes id_{V_2}}{\longrightarrow}
\Sim^2 (V_2) \otimes V_2
\stackrel{\mu_{2,1}}\rightarrow
\Sim^3 (V_2).
\end{multline}
In other words, writing $i_3$ for the composition of the maps in (\ref{i3def}),
  we have an exact sequence
\begin{equation}\label{i3}
0
\rightarrow
(\det V_1)^2 \otimes V_2
\stackrel{i_3}\longrightarrow
\Sim^3 (V_2)
\rightarrow
\A_6
\rightarrow
0.
\end{equation}

The $5$ ingredients are required to satisfy some open conditions,
just to ensure that what you cook is {\it eatable}. We need first to
give the recipe.

\subsection{The recipe}\label{recipe}
The conic bundle $\C$ comes from the first $4$ ingredients, and more
precisely is the image of the {\it relative $2$-Veronese}  $\PP(V_1)
\dashrightarrow \PP(V_2)$  given by the map $\sigma_2$ in the exact
sequence (\ref{conic}).

We give an {\it equation} defining $\C$ as a divisor in $\PP(V_2)$. A
conic bundle in a projective bundle $\PP(V)$ is given by an injection
of a line bundle to $\Sim^2 V$; in this case the map $\Sim^2(\sigma_2) \circ
\nu\colon (\det V_1)^2 \rightarrow \Sim^2V_2$.

\smallskip

Now we explain how to get $\Delta$ from $w$. The curve $\Delta$ is
{\it locally} (on $B$) the complete intersection of $\C$ with a relative cubic
in $\PP(V_2)$. In other words, a divisor in the linear system
associated to the restriction to $\C$ of the line bundle
$\hol_{\PP(V_2)}(3) \otimes \pi^* \LLL^{-1}$ for $\pi$ being the
projection on $B$, $\LLL$ a line bundle on $B$.

Why a map from a line bundle to the vector bundle $\A_6$ gives such a
divisor?
The {\it equation} of a divisor $\G \in |\hol_{\PP(V_2)}(3) \otimes \pi^*
\LLL^{-1}|$ is an injective  map $\LLL \hookrightarrow
\Sim^3V_2$. Intersecting it with $\C$ we do not obtain all divisor in
that linear system since in general they are not all complete
intersections of the form $\C \cap \G$. To get the complete linear
system we need to consider injections $\LLL \hookrightarrow
\A_6$ where $\A_6$ is the quotient of $\Sim^3 V_2$ by the subbundle
{\it of the relative cubics vanishing on $\C$}, that is exactly the
image of the map $i_3$ in the exact sequence (\ref{i3}).

\subsection{The open conditions}\label{open}
We need to impose that
\begin{itemize}
\item $\C$ has only Rational Double Points as singularities;
\item the curve $\Delta$ has only simple singularities,
   where ``simple'' means that the germ of double cover of $\C$ branched on it
   is either smooth or has a Rational Double Point.
\end{itemize}

\begin{df}\label{P}
The map $\sigma_2$ gives isomorphisms of the respective fibres of $\Sim^2 V_1$
and $V_2$ over points not in $\supp(\tau)$. On the points of
$\supp(\tau)$ it defines a rank $2$ matrix, whose image defines
a pencil of lines in the corresponding $\PP^2$, thus having a base
point.
  We denote by $\PPP$ the union of these (base) points. So $\PPP$ is
in natural bijection with $\supp(\tau)$.
\end{df}

\begin{oss}\label{DeltaP}
By theorem 4.7 of \cite{cp}, $\PPP \subset \Sing(\C)$ is the set of isolated
branch points of $\psi_2$, so
in particular $\Delta \cap \PPP = \emptyset$.
\end{oss}
\begin{oss}\label{easycond}
By remark 4.14 in \cite{cp}, if $\tau$ is a reduced divisor and every fibre of
$\C \rightarrow B$ is reduced (it is enough to check the preimages of points of
$\tau$, the other fibres being smooth) then the first open condition is
fulfilled. More precisely automatically $\Sing(\C)=\PPP$ and these points are
$A_1$ singularities of $\C$.

It follows that if moreover $\Delta$ is smooth and $\Delta \cap \PPP=\emptyset$
both open conditions are fulfilled and the relative canonical model of the
surface is smooth.
\end{oss}

\subsection{The dish}\label{invs}
What we get is a genus $2$ fibration $f\colon S \rightarrow B$ (the base is the
first ingredient) with $V_1 \cong f_*\omega_{S|B}$ and $V_2 \cong
f_*\omega_{S|B}^2$. The structure theorem says that any relatively minimal
genus $2$ fibration is obtained in this way.

Denoting by $b$ the genus of the base curve $B$
\begin{equation*}%\label{chi}
\chi(\hol_S)= \deg V_1 + (b-1)\ \ \ \
%\end{equation}
%\begin{equation}\label{Kquadro}
K^2_S=2\deg V_1 + 8(b-1)+\deg \tau
\end{equation*}

\section{The families}\label{8fam}

In this section we construct $8$ families of surfaces of general type with
$p_g=q=1$, $K^2=4$ and Albanese of genus $2$ using the recipe described in the
section \ref{struct}. We need then to give the ingredients, quintuples
$(B,V_1,\tau,\xi,w)$ with $B$ elliptic curve and (by \ref{invs}) $\deg V_1=1$,
$\deg \tau=2$.

\medskip

As first ingredient we take any elliptic curve $B$. For later
convenience we fix a group structure on $B$ and denote by $\eta_0=0$
its neutral element, and by $\eta_1$, $\eta_2$ and $\eta_3$ the
nontrivial 2-torsion points.

\medskip

The choice of the next $3$ ingredients for the $8$ families is summarized
in the tables \ref{sigdec} and \ref{sig}, which we are going to explain.

\medskip

As second ingredient, $V_1$, we need a vector bundle of rank $2$. $V_1$ can be
sum of line bundles (table \ref{sigdec}) or indecomposable (table \ref{sig}).

In the decomposable case we take  $V_1 \cong \hol_B(p) \oplus
\hol_B(0-p)$ where $p$ is a ${\mathfrak t}$-torsion point for some
${\mathfrak t} \in\{2,3,4,6\}$,
$V_2:=\hol_B(D_1) \oplus \hol_B(D_2) \oplus \hol_B(D_3)$ for
$D_1$, $D_2$ and $D_3$ suitable divisors on $B$. Since
  $$V_1 \cong \hol_B(p) \oplus \hol_B(0-p) \Rightarrow \Sim^2 V_1 \cong
\hol_B(2\cdot p) \oplus \hol_B(0) \oplus \hol_B(2\cdot 0-2\cdot p)$$
the splitting of the source and the target of $\sigma_2$ as sum of line
bundles allows to represent $\sigma_2$ by a $3 \times 3$ matrix whose
entries are global sections of line bundles over $B$. The table
\ref{sigdec} give $4$ families of choices of $\mathfrak t$, $D_1$, $D_2$,
$D_3$ and $\sigma_2$. The pair $(a_i,b_i)$ must be taken general in the sense of
\ref{open}, and we will later show that this open condition is
nonempty. The linear system on which $\tau$ varies depends on
the other data, and can be computed by (\ref{conic}): we wrote the result on
the last column.

\begin{table}
\caption{$\sigma_2\colon \Sim^2V_1 \rightarrow \bigoplus_{i=1}^3 
\hol_B(D_i)$ for $V_1 \cong \hol_B(p) \oplus \hol_B(0-p)$, $p$ 
$\mathfrak t$-torsion}\label{sigdec}
\renewcommand{\arraystretch}{1.3}
\begin{tabular}{|c||c|c|c|c|c|c|}
\hline
family&$\mathfrak t$&$D_1$&$D_2$&$D_3$&$\sigma_2$&$|\tau|$\\
\hline
\hline
$\M_{2,3}$&2&$2 \cdot 0$&$2\cdot0$&$0$&$\begin{pmatrix}
0&0&a_1\\
1&0&b_1\\
0&1&0
\end{pmatrix}$&$|2\cdot 0|$\\
$\M_{4,2}$&4&$2 \cdot 0$&$2\cdot p$&$0$&$\begin{pmatrix}
0&0&a_2\\
1&0&b_2\\
0&1&0
\end{pmatrix}$&$|2 \cdot p|$\\
$\M_{3,1}$&3&$0+ p$&$2 \cdot  p$&$0$&$\begin{pmatrix}
0&0&a_3\\
1&0&b_3\\
0&1&0
\end{pmatrix}$&$|2 \cdot 0|$\\
$\M_{6,1}$&6&$4\cdot p - 2 \cdot 0$&$2 \cdot p$&$0$&$\begin{pmatrix}
0&0&a_4\\
1&0&b_4\\
0&1&0
\end{pmatrix}$&$|2 \cdot 0|$\\
\hline
\end{tabular}
\end{table}

Otherwise we take $V_1$ to be the only indecomposable rank $2$ vector bundle
on $B$ with $\det V_1=\hol_B(0)$. By \cite{atiyah}, (as shown for the
analogous case $K_S^2=3$ in \cite{cp}) it follows that also in this
case $\Sim^2V_1$ is sum of line bundles, and more precisely
$$\Sim^2 V_1 \cong \hol_B(\eta_1) \oplus \hol_B(\eta_2) \oplus
\hol_B(\eta_3).$$ Therefore also in this case, writing
$V_2:=\hol_B(D_1) \oplus \hol_B(D_2) \oplus \hol_B(D_3)$ we
can represent $\sigma_2$ by a matrix. The table \ref{sig} give $4$
families of choices of $D_1$, $D_2$, $D_3$ and $\sigma_2$, and the resulting
$\tau$ (it moves in a pencil in all cases but the first);
in the last row $\sigma$ denotes a nontrivial $3$-torsion
point of $B$. $a_i,b_i,c_i,d_i$ are general in the sense of
\ref{open}.

\begin{table}
\caption{$\sigma_2\colon \Sim^2V_1 \rightarrow \bigoplus_{i=1}^3
     \hol_B(D_i)$ for $V_1$ indecomposable, $\det V_1 \cong
     \hol_B(0)$}\label{sig}
\renewcommand{\arraystretch}{1.3}
\begin{tabular}{|c||c|c|c|c|c|}
\hline
family&$D_1$&$D_2$&$D_3$&$\sigma_2$&$\tau$\\
\hline
\hline
$\M_{i,3}$&$2\cdot 0$&$2\cdot 0$&$\eta_3$&$\begin{pmatrix}
a_5&0&0\\
0&d_5&0\\
0&0&1
\end{pmatrix}$&$=\eta_1+\eta_2$\\
$\M_{i,2}$&$2\cdot 0$&$\eta_1+\eta_2$&$\eta_3$&$\begin{pmatrix}
a_6&b_6&0\\
c_6&d_6&0\\
0&0&1
\end{pmatrix}$&$\in|2\cdot 0|$\\
$\M_{i,2}'$&$2\cdot 0$&$0+\eta_1$&$\eta_3$&$\begin{pmatrix}
a_7&b_7&0\\
c_7&d_7&0\\
0&0&1
\end{pmatrix}$&$\in |0+\eta_2|$\\
%&&&&\\
$\M_{i,1}$&$0+\sigma$&$2\cdot\sigma$&$\eta_3$&$\begin{pmatrix}
a_8&b_8&0\\
c_8&d_8&0\\
0&0&1
\end{pmatrix}$&$\in |0+\eta_3|$\\
\hline
\end{tabular}
\end{table}

\medskip

Now that we have the first $4$ ingredients, we can construct the conic
bundle. The splitting of $V_2$ as sum of line bundles gives
relative coordinates on $\PP(V_2)$, by taking the injections
$y_i\colon \hol_B(D_i) \hookrightarrow V_2$. We
can use these coordinates to give equations of $\C \subset \PP(V_2)$.

\begin{lem}\label{eqconic}
The conic bundle $\C$ obtained by the ingredients given in a row of
the table \ref{sigdec} or \ref{sig} following the recipe
in \ref{recipe} has the equation given in the first column (and
corresponding row) of the table \ref{Delta}.
\end{lem}
\begin{proof}
As explained in \ref{recipe}, an equation of $\C$ is given by the map
$\Sim^2(\sigma_2) \circ \nu$, where $\nu$ is given in (\ref{nu}).

In the cases of table \ref{sigdec} $V_1$ is sum of two line bundles,
so we can use
the splitting to give two generators $x_0$, $x_1$ on each stalk.
When we write $\Sim^2 V_1 \cong
\hol_B(2\cdot p) \oplus \hol_B(0) \oplus \hol_B(2\cdot 0-2\cdot p)$
the first summand correspond to $x_0^2$, the second to $x_0x_1$, the
third to $x_1^2$. So by the expression of $\sigma_2$
$$
\left\{
\begin{array}{lll}
x_0^2 &\mapsto &y_2\\
x_0x_1& \mapsto& y_3\\
x_1^2 &\mapsto &a_iy_1+b_iy_2
\end{array}
\right .
$$
and the equation $(x_0)^2(x_1)^2=(x_0x_1)^2$ maps to
$y_2(a_iy_1+b_iy_2)=y_3^2$.

In the cases of table \ref{sig}, $V_1$ is indecomposable so we do not have
``global'' $x_0,x_1$. Anyway, as noticed in remark 6.13 of \cite{cp},
the map $\nu \colon \hol_B(2 \cdot 0) \rightarrow \Sim^2(\bigoplus
\hol_B(\eta_i))$ is given by a $6 \times 1$ matrix whose entries are
\begin{itemize}
\item[-] $0$ the three entries corresponding to the ``mixed terms''
   ($\hol_B(\eta_i+\eta_j)$ for $i \neq j$), since
   $i \neq j \Rightarrow \Hom(\hol_B(2\cdot 0),\hol_B(\eta_i+\eta_j))=0$
\item[-] isomorphisms the three entries corresponding to the pure powers
   ($\hol_B(\eta_i+\eta_i)$) since the Veronese image of $\PP^1$ in
   $\PP^2$ has rank $3$.
\end{itemize}
It follows that the equation of the relative Veronese embedding
$\PP(V_1) \hookrightarrow \PP(\Sim^2 V_1)$ is
$z_1^2+z_2^2+z_3^2=0$ for suitable choice of coordinates $z_i \colon
\hol_B(\eta_i) \hookrightarrow \Sim^2V_1$ on $\PP(\Sim^2
V_1)$. Composing with $\sigma_2$ we get the equations in the table.
\end{proof}

We still have to give the last ingredient.
Since in each case $\hol_B(2 \cdot \tau) \cong \hol_B(4 \cdot 0)$,
we have $(\det V_1 \otimes \hol_B(\tau))^2 \cong \hol_B(6 \cdot 0)$,
therefore $w$ is the class (modulo $\CC^*$) of a map $\hol_B(6 \cdot
0)\rightarrow \A_6$, where $\A_6$ is a quotient of $\Sim^3V_2$ as in
(\ref{i3}).

\begin{table}[ht]
\caption{$\C$ and $\Delta=\C\cap\G$}\label{Delta}
\renewcommand{\arraystretch}{1.3}
\begin{tabular}{|c||c|c|}
\hline
family&$\C$&$\G$\\
\hline
\hline
$\M_{2,3}$&$y_2(a_1y_1+b_1y_2)=y_3^2 $&$\sum_0^3k_iy_1^{3-i}y_2^i=0$\\
&&\\
$\M_{4,2}$&$y_2(a_2y_1+b_2y_2)=y_3^2 $&$y_1(k_0y_1^2+k_2y_2^2)=0$\\
&&\\
$\M_{3,1}$&$y_2(a_3y_1+b_3y_2)=y_3^2 $&$k_0y_1^3+k_3y_2^3=0$\\
&&\\
$\M_{6,1}$&$y_2(a_4y_1+b_4y_2)=y_3^2 $&$k_0y_1^3+k_3y_2^3=0$\\
&&\\
$\M_{i,3}$&$a_5^2y_1^2+d_5^2y_2^2+y_3^2=0$&$\sum_0^3k_iy_1^{3-i}y_2^i=0$\\
&&\\
$\M_{i,2}$&$(a_6y_1+c_6y_2)^2+(b_6y_1+d_6y_2)^2+y_3^2=0$&$y_1(k_0y_1^2+k_2y_2^2)=0$\\
&&\\
$\M_{i,2}'$&$(a_7y_1+c_7y_2)^2+(b_7y_1+d_7y_2)^2+y_3^2=0$&$y_1(k_0y_1^2+k_2y_2^2)=0$\\
&&\\
$\M_{i,1}$&$(a_8y_1+c_8y_2)^2+(b_8y_1+d_8y_2)^2+y_3^2=0$& 
$k_0y_1^3+k_3y_2^3=0$\\
\hline
\end{tabular}
\end{table}

We choose this map as composition of a general map $\overline{w}\colon \hol_B(6
\cdot 0) \rightarrow \Sim^3V_2$ with the projection to the quotient.  This
geometrically means that we take $\Delta=\C \cap \G$ for a relative cubic
$\G\subset \PP(V_2)$ whose equation is given by $\overline{w}$. Since
$\Sim^3V_2$ is sum of line bundles whose maximal degree is $6$, the nonzero
entries of $\bar{w}$ are constants and correspond to the summands of the target
isomorphic to $\hol_B(6\cdot 0)$. In the table \ref{Delta} we give the exact
equation of $\G$ in each case. The parameters $k_i\in \CC$ must be taken
general in the sense of \ref{ingr}, requiring that $\Delta$ has only simple
singularities.

\smallskip

\begin{prop}\label{4famiglie!}
Cooking the ingredients given above ($B$ general elliptic curve, $V_1$, $\tau$,
$\xi$ given by a row of the table \ref{sigdec} or \ref{sig}, $w$ by the
corresponding row in the table \ref{Delta}) following the recipe \ref{recipe},
one finds $8$ unirational families of minimal surfaces of general type with
$p_g=q=1$, $K^2=4$, Albanese morphism $\alpha$ with fibres of genus $2$ and
$\alpha_*\omega^2_S$ sum of line bundles. The general element in each family
has ample canonical class.
\end{prop}

\begin{proof}
By the recipe (\ref{recipe}) and remark \ref{easycond}, if we show
that all these families of ingredients contain one element such that
\begin{itemize}
\item[on $\tau$:] $\coker \sigma_2 \cong \hol_\tau$ for $\tau$ reduced divisor;
\item[on $\C$:] all fibres of $\C \rightarrow B$ are reduced conics;
\item[on $\Delta$:] $\Delta$ is smooth and $\Delta \cap \PPP =\emptyset$.
\end{itemize}
then all these examples give families of genus $2$ fibrations $f\colon
S \rightarrow B$ with (by \ref{invs}) $K_S^2=4$ and $\chi(\hol_S)=1$
with smooth relative canonical model. Since $B$ has genus $1$, $q(S) \geq 1$,
so $p_g=q=1$. By the universal property of the Albanese morphism $\alpha=f$,
and therefore $\alpha_*\omega_S \cong V_1$, $\alpha_*\omega^2_S \cong
\hol(D_1) \oplus\hol(D_2) \oplus \hol(\eta_3)$.

So we only need to find an element in each family satisfying the three
condition. Since all conditions are open and each family irreducible, it is
enough to show that each condition (separately) is fulfilled by some choice of
the parameters. This is easy, we sketch a way to do it.
\begin{itemize}
\item[On $\tau:$] we need to choose the entries of the matrix of $\sigma_2$ so
that the determinant is not a perfect square.
\item[On $\C:$] a conic of the form $y_3^2=q(y_1,y_2)$ is a double line if and
only if $q=0$. By the equation of $\C$ in the table \ref{Delta} we see that in
the first $4$ cases it is enough to choose $a_i,b_i$ without common zeroes,
whereas in the last $4$ cases it is enough $\det \sigma_2 \neq 0$.
\item[On $\Delta$:] in $5$ cases the linear system $|\G|$ has fixed locus
$\{y_1=y_2=0\}$ which do not intersect $\C$. So $|\Delta|$ is free and
therefore we can conclude by Bertini. In the remaining cases $\M_{4,2}$,
$\M_{i,2}$ and $\M_{i,2}'$ the fixed part of $|\Delta|$ is $\{y_1=0\} \cap \C$
and the general element of the movable part of $|\Delta|$ do not intersect the
fixed part. So we only need to check $\{y_1=0\} \cap \C$ smooth and not
containing $\PPP$. For $\M_{4,2}$ if we take $b_2\neq 0$ we get smoothness,
and the other condition comes automatically since $\PPP \subset
\{y_2=y_3=0\}$. In the other two cases $c_i^2+d_i^2$ square free gives the
smoothness, and from ({\it e.g.}) $a_ib_ic_id_i\neq 0$ follows 
$\{y_1=0\} \cap \C
\not\supset \PPP$.
\end{itemize}
\end{proof}

We end the section by explaining the choice of the indices of the name
of each family.

The first index remembers us which $V_1$ we have chosen: $i$ stands for ``$V_1$
indecomposable'', a number $t$ means ``$V_1$ has a $t$-torsion bundle as direct
summand''.

The second index gives the number of connected components of the curve
$\Delta$ for a surface in the family. Let us show this decomposition.

The equation of $\G$ is homogeneous of degree $3$ in two variables
(with constant coefficients), so we can formally decompose it as
product of three linear factors. When $D_1=D_2$ ($\M_{2,3}$ and $\M_{i,3}$)
each factor gives a map of a line bundle ($\hol_B(2\cdot 0)$) to
$V_2$, so a relative hyperplane of $\PP(V_2)$: these three relative
hyperplanes cut on $\C$ three components of $\Delta$ that pairwise
they do not intersect.

When $\hol_B(D_1) \not\cong\hol_B(D_2)$ a factor $cy_1+c'y_2$
determines a relative hyperplane only if $cc'=0$. In the cases $\M_{4,2}$,
$\M_{i,2}$, $\M_{i,2}'$ one can then decompose $\Delta$ as union of its fixed
part $\{y_1=0\}$ and its movable part.

\section{Direct image of the canonical sheaf decomposable}\label{V1dec}
In this section we prove the following

\begin{prop}\label{v1dec}
All minimal surfaces of general type $S$ with $K_S^2=4$, $p_g=q=1$ 
such that the general fibre of the Albanese morphism $\alpha$
has genus $2$ and $\alpha_* \omega_S$, $\alpha_* \omega_S^2$ are
direct sum of line bundles belong to $\M_{2,3}$, $\M_{4,2}$, $\M_{3,1}$
or $\M_{6,1}$.
\end{prop}

By the structure theorem of genus $2$ fibrations, we need to classify
$5$-tuples $(B,V_1,\tau,\xi,w)$ with $B$ elliptic curve,
$\deg V_1=1$, $\deg \tau=2$ such that $V_1$ and $V_2$ are
sum of line bundles.

Since $h^0(V_1)=h^0(\omega_S)=p_g$ we can assume up to translations
$V_1\cong \hol_B(p) \oplus \hol_B(0-p)$ for some $p \neq 0$. We write
$V_2 = \hol_B(D_1) \oplus \hol_B(D_2) \oplus \hol_B(D_3)$, with $D_i$
divisors of degree $d_i$, $d_3 \leq d_2 \leq d_1$. We consider
relative coordinates in $V_1$ and $V_2$ as follows: $x_i$ correspond to
the summand of degree $i$ in $V_1$, $y_j$ correspond to the summand
$\hol_B(D_j)$ in $V_2$.

\begin{lem}\label{221}
$d_1=d_2=2$, $d_3=1$.
\end{lem}
\begin{proof}
By the exact sequence (\ref{conic}), since $\Sim^2V_1$ is direct sum
of three line bundles of respective degrees $0$, $1$ and $2$,
$d_1+d_2+d_3 =5$, $d_i \geq 3-i$.

Since $d_3 \leq d_2 \leq d_1$ to show $d_3=1$ we assume by contradiction
$d_3=0$. Then the summands of positive degree in $\Sim^2 V_1$ map trivially on
$\hol_B(D_3)$. In other words $\sigma_2(x_1^2), \sigma_2(x_0x_1) \in
\Span(y_1,y_2)$. In particular, the equation of $\C$ being
$\sigma_2(x_0^2)\sigma_2(x_1^2)=\sigma_2(x_0x_1)^2$, the section
$s:=\{y_1=y_2=0\}$ is contained in $\C$.

We consider $s$ as Weil divisor in $\C$. Note that $\C$ has only canonical
singularities, so $s$ is $\QQ$-Cartier, and the self-intersection number $s^2$
is well defined, as the numbers $s\cdot D \in \ZZ$ for any Cartier divisor $D$
on $\C$, including $K_{\C}$.

We denote by $H$ the numerical class of a divisor in $\hol_{\PP(V_2)}(1)$, by
$F$ the class of a fiber of the map $\PP(V_2) \rightarrow B$. Then $s$, as a
cycle in $\PP(V_2)$, has numerical class $(H-d_1F)(H-d_2F)=H^2-5HF$. $\Delta$
is Cartier on $\C$, and the corresponding line bundle is the restriction to
$\C$ of a line bundle in $\PP(V_2)$ whose numerical class is $3H-6F$ (since
$\deg (\det V_1 \otimes \hol_B(\tau))^2=6$). It follows $\Delta \cdot
s=(3H-6F)(H^2-5HF)=-6<0$. Being $s$ irreducible, $s<\Delta$.

Consider now a minimal resolution of the singularities $\rho\colon \tilde{\C}
\rightarrow \C$ and let $\tilde{s}$ be the strict transform of $s$. Then
$\tilde{s}$ is a smooth elliptic curve and $\tilde{s}=\rho^*s-e$ for some
exceptional $\QQ$-divisor $e$, so $s^2+K_{\C}s\geq
s^2+e^2+K_{\C}s=\tilde{s}^2+K_{\tilde{C}}\tilde{s}=0$. Since the class of $\C$
is $2H-2F$ it follows $-s^2 \leq K_{\C}s = (-H+3F)(H^2-5HF)=3$ and
$(\Delta-s)s=-s^2-6\leq -3<0$. It follows so $2s < \Delta$, contradicting
\ref{open}.

Then $d_3=1$ and to conclude we can assume by contradiction $d_2=1$, then
$\sigma_2(x_1)^2\in \Span(y_1)$. It follows that the equation of $\C$ is a
square modulo $y_1$. In other words the relative hyperplane $\{y_1=0\}$ cut
$2\cdot s'$ where $s$ is a section of the map $\PP(V_2) \rightarrow B$. The
class of $s'$ is $H^2-4HF$: repeating the above argument we find $\Delta \cdot
s' =-3$, $(\Delta - s')\cdot s' \leq -1 \Rightarrow 2s'<\Delta$, the same
contradiction as above.
\end{proof}

\begin{lem}\label{x0x1y3}
$\sigma_2(x_0x_1) \not\in \Span (y_1,y_2)$.
\end{lem}

\begin{proof}
Since $\sigma_2(x_1^2) \in \Span (y_1,y_2)$, if also $\sigma_2(x_0x_1)
\in \Span (y_1,y_2)$, then  the section  $s:=\{y_1=y_2=0\}$ is
contained in $\C$. The numerical class of $s$ is $H^2-4HF$ so (as in
the previous proof) $\Delta \cdot s =-3$, $(\Delta - s)\cdot s \leq -1
\Rightarrow 2s<\Delta$, a contradiction.
\end{proof}

\begin{oss}\label{D3=0}
The lemma \ref{x0x1y3} says that the composition of $\sigma_2$ with the
projection onto the summand $\hol_B(D_3)$ is different from zero. Since any
nonzero morphism between line bundles of the same degree is an isomorphism, it
follows $\hol_B(D_3) \cong \hol_B(0)$.
\end{oss}

\begin{lem}\label{splitV1dec}
The exact sequence (\ref{i3}) splits.
\end{lem}
\begin{proof}
By the lemma \ref{x0x1y3} and remark \ref{D3=0} the coefficient of the
term $y_3^2$ in the relative conic $\sigma_2(x_0^2)\sigma_2(x_1)^2 -
\sigma_2(x_0x_1)^2$ defining $\C$ is a nonzero constant. Then each
relative conic can be uniquely decomposed as a sum of a multiple of this
equation with an equation where the multiples of $y_3^2$  ($y_1y_3^2,
y_2y_3^2, y_3^3$) do not appear.

Since the multiples of the equation of $\C$ define exactly the image
of $i_3$, this means that the restriction of the projection $\Sim^3V_2
\rightarrow \A_6$ to $\Sim^3(\hol_B(D_1) \oplus \hol_B(D_2)) \oplus
(\Sim^2(\hol_B(D_1) \oplus \hol_B(D_2))\otimes \hol_B(D_3))$ is an
isomorphism. Its inverse splits the exact sequence (\ref{i3}).
\end{proof}

In particular every morphism to $\A_6$ lift to a morphism to
$\Sim^3V_2$, and therefore the last ``ingredient'' $w$ comes from a
map $\bar{w}\colon (\det V_1 \otimes \hol_B(\tau))^2 \rightarrow
\Sim^3V_2$. It follows

\begin{cor}\label{t-tors}
$\T:=\hol_B(D_1-D_2)$ is a $\mathfrak{t}$-torsion bundle for some 
$\mathfrak{t}\in
\{1,2,3\}$, and up to exchange $D_1$ and $D_2$, $\hol_B(0+\tau)^2\cong
\hol_B(D_1)^3$.
\end{cor}
\begin{proof}
The source of $\bar{w}$ is the line bundle $\hol_B(0+\tau)^2$ of
degree $6$. Since $\Sim^3V_2$ is sum of line bundles of degree at
most $6$, its image is contained in the sum of those having exactly degree $6$,
$\Sim^3(\hol_B(D_1) \oplus \hol_B(D_2))$, and more precisely in those
summands isomorphic to $\hol_B(0+\tau)^2$.

So $\Delta =\C \cap \G$ with $\G = \sum k_i y_1^{3-i} y_2^i$ where
$k_i$ are constant that can be different form $0$ only when
$\hol_B((3-i)D_1+iD_2) \cong \hol_B(0+\tau)^2$.
The claim follows since $\Delta$ is reduced, and then at least two $k_i$'s
are different from $0$.
\end{proof}

\begin{proof}[Proof of proposition \ref{v1dec}]
By remark \ref{D3=0} and corollary \ref{t-tors} $V_2\cong \T(D_2)
\oplus \hol_B(D_2) \oplus \hol_B(0)$ for some $\mathfrak{t}$-torsion line
bundle $\T$, $\mathfrak{t}\in \{1,2,3\}$. Moreover, by exact sequence
(\ref{conic}) and corollary \ref{t-tors}
$$
\T(2\cdot D_2+0) \cong \hol_B(3\cdot 0 +\tau)\ \ \ \ \
\hol_B(2\cdot 0+2 \cdot \tau)\cong \T^3(3\cdot D_2).
$$
equivalently
\begin{equation}\label{formulette}
\hol_B(D_2) \cong \T(2 \cdot 0)\ \ \ \ \
\hol_B(\tau) \cong \T^3(2\cdot 0)
\end{equation}

Moreover, by the injectivity of $\sigma_2$, $2p$ must be linearly
equivalent to $D_1$ or $D_2$, {\it i.e.}
\begin{equation}\label{2p=Di}
\hol_B(2\cdot p) \cong \T(2 \cdot 0)\ \
\text{ or }\ \ \hol_B(2p) \cong \T^2(2 \cdot 0)
\end{equation}

\begin{itemize}
\item[If $\mathfrak{t}=1$:] $\T\cong \hol_B$ and the two alternatives in
   (\ref{2p=Di}) are identical: $\hol_B(2\cdot p) \cong \hol_B(2\cdot 0)$. Since
   $p \neq 0$, $p$ is a $2$-torsion point. We can choose coordinates in $V_2$
   such that $y_2=\sigma_2(x_1^2)$ and (by lemma \ref{x0x1y3})
   $y_3=\sigma_2(x_0x_1)$. We can also assume $\sigma_2(x_0^2) \in \Span
   (y_1,y_2)$ by changing the coordinates $(x_0, x_1$): we have found the
   family $\M_{2,3}$.
\item[If $\mathfrak{t}=2$:] If $\hol_B(2\cdot p) \cong \T(2 \cdot 0)$, $p$ is a
   $4$-torsion point. Changing coordinates in $V_1$ and $V_2$ as above
   we find the family $\M_{4,2}$. \\ Else
   $\hol_B(2\cdot p) \cong \T^2(2 \cdot 0)$. In this case ($\hol_B(D_2)
   \not\cong \hol_B(2\cdot p)$) $\sigma_2(x_1^2) \in \Span(y_1)$,
   therefore (see definition \ref{P}) $\PPP \subset \{y_1=0\}$. On the
   other hand $\G=\{y_1(k_0y_1^2+k_2y_2^2)=0\}$, so the fixed part of
   $|\Delta|$ contains $\PPP$, contradicting remark \ref{DeltaP}: this
   case do not occur.
\item[If $\mathfrak{t}=3$:] If $\hol_B(2\cdot p) \cong \T(2 \cdot 0)$, $p$ is
   either a $3$-torsion point or a $6$-torsion point. Changing coordinates as
   above we find respectively the families $\M_{3,1}$ and $\M_{6,1}$. The other
   case $\hol_B(2\cdot p) \cong \T^2(2 \cdot 0)$ gives the same families (with
   $D_1$ and $D_2$ exchanged).
\end{itemize}
\end{proof}

\section{Direct image of the canonical sheaf indecomposable}\label{V1undec}
In this section we prove the following

\begin{prop}\label{v1undec}
All minimal surfaces of general type $S$ with $K_S^2=4$, $p_g=q=1$
such that the general fibre of the Albanese morphism $\alpha$ has
genus $2$, $\alpha_* \omega_S$ is an indecomposable vector bundle and
$\alpha_* \omega_S^2$ is a direct sum of line bundles belong to
$\M_{i,3}$, $\M_{i,2}$, $\M_{i,2}'$ or $\M_{i,1}$.
\end{prop}

We need to classify $5$-tuples $(B,V_1,\tau,\xi,w)$ with $B$ elliptic curve,
$V_1$ indecomposable of degree $1$, $\deg \tau=2$ such that $V_2$ is
sum of three line bundles.

$B$ can be any elliptic curve and by Atiyah's classification of the
vector bundles on an elliptic curves \cite{atiyah}, we can assume (up to
translations) $V_1=E_0(2,1)$, that is the only indecomposable vector
bundle over $B$ whose determinant is $\hol_B(0)$.

>From Atiyah's results follows $\Sim^2 (V_1) \cong \hol_B(\eta_1) \oplus
\hol_B(\eta_2) \oplus \hol_B(\eta_3)$. As in the previous case we write $V_2 =
\hol_B(D_1) \oplus \hol_B(D_2) \oplus \hol_B(D_3)$, with $D_i$ divisors of
degree $d_i$, $d_3 \leq d_2 \leq d_1$.

\begin{oss}\label{conicundec}
As shown in the proof of lemma \ref{eqconic}, in this case the relative
$2$-Veronese $\PP(V_1) \hookrightarrow \PP(\Sim^2V_1)$ has equation
$z_1^2+z_2^2+z_3^2=0$ for a suitable choice of coordinates $z_i \colon
\hol_B(\eta_i) \hookrightarrow \Sim^2V_1$.

It follows that, in these coordinates, $\C$ is defined by the polynomial
$\sum_{i=1}^3 \sigma_2(z_i)^2$.
\end{oss}

\begin{lem}\label{splits}
We can assume $D_3=\eta_3$, and we can choose coordinates in $V_2$ so that
$\sigma_2(z_3)=y_3$. Moreover the exact sequence (\ref{i3}) splits.
\end{lem}

\begin{proof}
Since $\sum d_i=5$ and by the injectivity of $\sigma_2$, $\forall i\
d_i\geq 1$, $d_3=1$. The injectivity of $\sigma_2$ forces now one
of the induced maps $\hol_B(\eta_i) \rightarrow \hol_B(D_3)$ to be an
isomorphism and then (renaming the torsion points) we have
$D_3=\eta_3$. Changing coordinates in $V_2$ we can assume
$\sigma_2(\hol_B(\eta_3))=\hol_B(D_3)$.

By remark \ref{conicundec} the coefficient of the term $y_3^2$ in the
equation of $\C$ is a nonzero constant and we can conclude as in the proof
of lemma \ref{splitV1dec}.
\end{proof}

\begin{lem}\label{dec}
$d_1=d_2=2$.
\end{lem}

\begin{proof}
We assume by contradiction $d_2=1$, $d_1=3$. By lemma \ref{splits} the
curve $\Delta$ is a complete intersection $\G \cap \C$ for a relative
cubic $\G$ defined by an immersion $\bar{w}$ of a line bundle of degree $6$ to
$\Sim^3V_2$.

The image of $\overline{w}$ is then contained in $\hol_B(D_1)^2
\otimes V_2$ since all other summands have degree strictly smaller
than $6$. In other words the equation of $\G$ is divisible by
$y_1^2$. In particular $\Delta$ contains $\{y_1=0\} \cap \C$ with
multiplicity $2$, contradicting \ref{open}.
\end{proof}

It follows, as in the previous case

\begin{cor}\label{t-torsundec}
$\T:=\hol_B(D_1-D_2)$ is a $\mathfrak{t}$-torsion bundle for some
$\mathfrak{t}\in\{1,2,3\}$ and, up to exchange $D_1$ and $D_2$,
$\hol_B(0+\tau)\cong \hol_B(D_1)^3$.
\end{cor}

\begin{proof}
Identical to the proof of the analogous corollary \ref{t-tors}.
\end{proof}

\begin{proof}[Proof of proposition \ref{v1undec}]
By lemma \ref{splits} and corollary \ref{t-torsundec}, $V_2 \cong
\T(D_2) \oplus \hol_B(D_2) \oplus \hol_B(\eta_3)$, and, by the exact
sequence (\ref{conic}) and corollary \ref{t-torsundec}
$$
\T(2\cdot D_2+\eta_3)\cong \hol_B(3\cdot 0 +\tau)\ \ \ \ \
\hol_B(2\cdot 0+2\cdot \tau) \cong \T^3(3 \cdot D_2),
$$
equivalently
\begin{equation}\label{formuletteundec}
\hol_B(D_2) \cong \T(2\cdot 0) \ \ \ \ \
\hol_B(\tau) \cong \T^3(0 +\eta_3).
\end{equation}
Recall that by lemma \ref{t-torsundec} we can choose
$y_3=\sigma_2(z_3)$ and since $d_1=d_2=2$,  $\sigma_2(z_1)$,
$\sigma_2(z_2) \in \Span(y_1,y_2)$. In other words the matrix of
$\sigma_2$ is as the matrices in the last three rows of
table \ref{sig}.
\begin{itemize}
\item[If $\mathfrak{t}=1$:] $\T\cong \hol_B$, $\hol_B(D_1) \cong \hol_B(D_2)
   \cong \hol_B(2\cdot 0)$ and $\hol_B(\tau)\cong \hol_B(0+\eta_3)$. In
   fact, since $D_1=D_2$ we can change coordinates in $V_2$ to add
   to one of the first two rows any multiple of the other and diagonalize
   the matrix: this is the family $\M_{i,3}$. Note that $\tau=\eta_1+\eta_2$
   cannot move.
\item[If $\mathfrak{t}=2$:] then either $\T \cong \hol_B(\eta_3)$ or we can
   rename $\eta_1$ and $\eta_2$ to get $\T \cong \hol_B(\eta_1)$. This
   gives respectively the families $\M_{i,2}$ and $\M_{i,2}'$.
\item[If $\mathfrak{t}=3$:] Then $\T\cong \hol_B(0-\sigma)$ for some 
$3$-torsion
   point $\sigma$. This is the family $\M_{i,1}$.
\end{itemize}
\end{proof}

\section{Moduli}\label{def}

In this section we consider the scheme $\M$ in theorem \ref{main},
subscheme of the moduli space of the minimal surfaces of general type
given by the surfaces with $p_g=q=1$, $K^2=4$ whose Albanese fibration
$\alpha$ has general fibre a genus $2$ curve and such that
$\alpha_*\omega_S^2$ is sum of line bundles.

We have constructed $8$ unirational families of such surfaces in
proposition \ref{4famiglie!}, labeled $\M_{2,3}$, $\M_{4,2}$,
$\M_{3,1}$, $\M_{6,1}$, $\M_{i,3}$, $\M_{i,2}$, $\M_{i,2}'$ and
$\M_{i,1}$. Their parameter spaces have a natural map to $\M$.

\begin{oss}\label{comp}
$\M$ has $8$ connected components, that with a natural abuse of
notation we will denote by $\M_{2,3}$, $\M_{4,2}$, $\M_{3,1}$,
$\M_{6,1}$, $\M_{i,3}$, $\M_{i,2}$, $\M_{i,2}'$ and $\M_{i,1}$. Each
component is the image of the parameter space of the namesake family,
in particular is unirational.
\end{oss}

\begin{proof}
The map from the parameter space of our families to $\M$ is surjective
by propositions \ref{v1dec} and \ref{v1undec}.

There are many way to show that the closure of the images of two of
these parameter spaces do not intersect. For example, since the number of
direct summands of $V_1$ is a topological invariant by \cite{cc1},
$$(\overline{\M_{2,3} \cup \M_{4,2} \cup \M_{3,1} \cup \M_{6,1}}) \cap
(\overline{\M_{i,3} \cup \M_{i,2} \cup \M_{i,2}' \cup \M_{i,1}})=\emptyset.$$

The closure of two of the first $4$ families cannot intersect
because the degree $0$ summand of $V_1$ is in all cases a torsion line
bundle but with different torsion order. 
To show $\overline{\M_{i,2}}\cap \overline{\M_{i,2}'}=\emptyset$ we apply
the same argument to $(\det V_1)^2 \otimes \hol_B(-\tau)$. Finally the
same argument applied to $\hol_B(D_1-D_2)$ shows that also the
closures of the remaining pairs of families do not intersect. 

\end{proof}

\begin{prop}\label{dimension}
$\dim \M_{2,3}=5$. All other components of $\M$ have dimension $4$.
\end{prop}
\begin{proof}
The natural way to compute the dimension of each component is
computing the dimension of the corresponding parameter space, and then
subtract to the result the dimension of the general fibre of the map
into $\M$. These fibres correspond to orbits for the action of certain
automorphism groups.

$\Aut V_1$ and $\Aut V_2$ do not act on our data, since in the tables
\ref{sigdec} and \ref{sig} we require the matrix of $\sigma_2$ to have
special form. But in fact in all cases this ``special'' form is the
form of a general morphism $\Sim^2V_1 \rightarrow V_2$ in suitable
coordinates (for $V_1$ and $V_2$). It is then equivalent (but easier
to compute) to consider $\sigma_2$ general in $\Hom(\Sim^2 V_1, V_2)$
and act on it with the full group $\Aut V_1 \times \Aut V_2$.

Are there other automorphisms to consider? We can forget the action of $\Aut B$
since we have fixed a point of $B$ by choosing $\det V_1 \cong \hol_B(0)$, so
only a finite subgroup of $\Aut B$ act on our data, and quotienting by it do
not affect the dimension. The other automorphism to consider is (since we are
interested in $\Delta$ and not in its equation) ``multiply the equation of $\G$
by a constant leaving the other data fixed''. If you prefer, that's the action
of the automorphisms of the line bundle $(\det V_1 \otimes
\hol_B(\tau))^2$. Anyway, multiplying $V_1$ by $\lambda$ and $V_2$ by
$\lambda^2$ do not change $\sigma_2$ but multiply the equation of $\G$ by
$\lambda^{-6}$: this shows that we can restrict to consider the action of $\Aut
V_1 \times \Aut V_2$.

We leave to the reader the check that the subgroup of $\Aut V_1 \times \Aut
V_2$ fixing our data is finite. It follows (the moduli space of elliptic curves
has dimension $1$) that the dimension of each family is
$$
1+h+\delta-\alpha_1-\alpha_2
$$
where $h$, $\delta$, $\alpha_i$ are respectively the dimensions of
$\Hom(\Sim^2V_1,V_2)$, $\Hom((\det V_1 \otimes \hol_B(\tau))^2,\Sim^3V_2)$ and
$\Aut V_i$.

Now the computation is easy:
$$\begin{array}{lrrrrrrrrrrr}
\dim \M_{2,3}= &1&+&10&+&4&-&3&-&7&=&5\\
\dim \M_{4,2}= &1&+& 9&+&2&-&3&-&5&=&4\\
\dim \M_{3,1}= &1&+& 9&+&2&-&3&-&5&=&4\\
\dim \M_{6,1}= &1&+& 9&+&2&-&3&-&5&=&4\\
\dim \M_{i,3}= &1&+& 7&+&4&-&1&-&7&=&4\\
\dim \M_{i,2}= &1&+& 7&+&2&-&1&-&5&=&4\\
\dim \M_{i,2}'=&1&+& 7&+&2&-&1&-&5&=&4\\
\dim \M_{i,1}= &1&+& 7&+&2&-&1&-&5&=&4
\end{array}$$
\end{proof}

\begin{prop}\label{moduli}
All connected components of $\M$ are irreducible components of the
moduli space of minimal surfaces of general type.
\end{prop}
\begin{proof}
We need to show that for the general surface in each component, $h^1(\T_S)$ is
not greater than the dimension of the family, say $d$. By proposition
\ref{dimension}, $d\in \{4,5\}$ and more precisely $d=5$ only for the family
$\M_{2,3}$.

Equivalently (by Serre duality and since $h^0(\T_S)=0$ for a surface of general
type) we can show $h^0(\Omega_S^1 \otimes \omega_S)= 2K_S^2-10\chi(\hol_S) +
h^1(\T_S)\leq d-2$.

For a fibration $f \colon S \rightarrow B$, we denote by $\Crit(f)
\subset S$ the scheme of its critical points, $\D \subset \Crit(f)$ its
divisorial part. By definition $\D$ is supported on the nonreduced
components of the singular fibres.

Then (cf. \cite{FabDispense} lect. 9) computing kernel and cokernel of the
natural map $\xi'\colon \Omega^1_S \rightarrow \omega_{S|B}$ locally defined by
$\xi'(\eta)=(\eta \wedge dt) \otimes (dt)^{-1}$ (for $t$ a local parameter on
$B$) one finds an exact sequence
\begin{equation}\label{longexactfab}
0 \rightarrow \hol_S (f^* \omega_B+\D)
\rightarrow \Omega^1_S
\rightarrow
\omega_{S|B}
\rightarrow
\hol_{\Crit(f)}(\omega_{S|B})
\rightarrow
0
\end{equation}

By the proof of proposition \ref{4famiglie!}, the Albanese fibration
$\alpha$ of a general element $S$ in each of our families factors as
composition of
\begin{itemize}
\item a conic bundle $\C \rightarrow B$ with two singular fibres, both
   reduced,  with $\Sing(\C)$ consisting in two nodes, at the vertices
   of the two singular fibres;
\item a finite double cover $S \rightarrow \C$ branched on the two nodes of
$\C$ and on a smooth curve $\Delta$ not passing through the nodes.
\end{itemize}
It follows that each component of each fibre of $\alpha$ is reduced, so
$\D=\emptyset$. Since $\omega_B=\hol_B$ twisting the exact sequence
(\ref{longexactfab}) by $\omega_S$ we get the exact sequence
$$
0 \rightarrow \omega_S
\rightarrow \Omega^1_S \otimes \omega_S
\rightarrow
\omega_S^2
\rightarrow
\hol_{\Crit(\alpha)}(\omega_{S}^2)
\rightarrow
0$$
Since $p_g=1$ the required inequality $h^0(\Omega_S^1 \otimes \omega_S)\leq
d-2$ follows if we show $\dim \ker \left(  H^0(\omega_S^2) \rightarrow
H^0(\hol_{\Crit(\alpha)}(\omega_{S}^2))\right)=d-3$. In other words we must
show that
\begin{itemize}
\item[1)] the set of bicanonical curves containing the $0$-dimensional
   scheme $\Crit(\alpha)$ of the general surface in $\M_{2,3}$ is a pencil;
\item[2)] the general surface in each of the other families has only one
   bicanonical curve containing $\Crit(\alpha)$.
\end{itemize}

We study the bicanonical system of $S$. The involution on a surface induced by
a genus $2$ fibration (acting as the hyperelliptic involution on any fibre)
acts on $H^0(2K_S)$ as the identity. In our cases, at least for a general
surface as above (the relative canonical model is smooth and minimal), the
quotient by this involution is $\C$. So the bicanonical system of $S$ is the
pull-back of a linear system on $\C$, more precisely ($\omega_S=\omega_{S|B}$)
the restriction of $|\hol_{\PP(V_2)}(1)|$.

We study the critical points of $\alpha$. Since $\C$ has only reduced fibres
the critical points of $\alpha$ must be fixed points for the involution on
$S$. The isolated fixed points are the preimages of the two nodes of $\C$, and
they are critical for $\alpha$ (in suitable local coordinates
$\alpha(x,y)=xy$). The other critical points of $\alpha$ lies on the divisorial
fixed locus of the involution, where the involution has the form $(x,y) \mapsto
(x,-y)$: they are critical for $\alpha$ if and only if $\frac{\partial
\alpha}{\partial x}=0$. In other words we need their image on $\C$ to be a
ramification point for the map $\Delta \rightarrow B$.

So we need to compute the dimension of the subsystem of $|\hol_{\PP(V_2)}(1)|$
containing the nodes of $\C$ and the critical points of the map $\Delta
\rightarrow B$. Note that by the local computation above this is true
schematically: we need $H$ to contain the zero dimensional scheme $\Sing(\C)
\cup \Crit(\Delta \rightarrow B)$.

In all cases (see table \ref{Delta}) $\C=\{q(y_1,y_2)+y_3^2=0\}$: in
particular the nodes of $\C$ lie in $\{y_3=0\}$. Moreover $\Delta=\C
\cap \G$ for $\G=\{G(y_1,y_2)=0\}$. $\Crit(\Delta \rightarrow B)$ is defined by
$$
\rank
\begin{pmatrix}
\frac{\partial q}{\partial y_1}&
\frac{\partial q}{\partial y_2}&
2y_3\\
\frac{\partial G}{\partial y_1}&
\frac{\partial G}{\partial y_2}&
0\\
\end{pmatrix}\leq 1
$$
therefore (being $q$ and $G$ homogeneous in the $y_i$'s) $\Crit(\Delta
\rightarrow B)=\Delta \cap \{y_3=0\}$.

We have shown that $\left( \Sing(\C) \cup \Crit(\Delta \rightarrow B)
\right) \subset \{y_3=0\}$.  First consequence is that any relative
hyperplane of the form $\{fy_3=0\}$ contains the nodes of $\C$ and
$\Crit(\Delta \rightarrow B)$.

Choosing $f \in H^0(\hol_B(D_3))$, $\hol_B(D_3)$ being the direct
summand of $V_2$ given by the coordinate $y_3$, we find a curve whose
pull-back is a bicanonical curve through $\Crit(\alpha)$.
Note that $\deg D_3=1$ so in all cases we have found exactly one
bicanonical curve through $\Crit(\alpha)$.

If there are further bicanonical curves through
$\Crit(\alpha)$, then in the corresponding system of relative
hyperplanes in $\PP(V_2)$ there is an element $H$ not containing
$\{y_3=0\}$ and $H \cap \C \cap \{y_3=0\}$ contains the
$0-$dimensional scheme $\Delta \cap \{y_3=0\}$.
If $H \cap \C \cap \{y_3=0\}$ is also $0$-dimensional, then
by intersection computation both $H \cap \C \cap \{y_3=0\}$ and
$\Delta \cap \{y_3=0\}$ have length $6$, so they must be equal, a contradiction
since $\Sing \C \subset H \cap \C \cap \{y_3=0\}$ but $\Sing(\C)
\not\subset \Delta$. Therefore, if there are further bicanonical curves through
$\Crit(\alpha)$, then $H \cap \C \cap \{y_3=0\}$ contains a curve.

To conclude the proof we must now argue differently according to the family.
\begin{itemize}
\item[($\M_{i,1}$, $\M_{i,2}'$, $\M_{i,3}$)] We set $b_5:=c_5:=0$ to
   treat these cases together. If $a_j,b_j,c_j,d_j$ have no common zeroes, $\C
   \cap \{y_3=0\}$ has a finite map of degree $2$ onto $B$ and then, if it is
   reducible, its components are cut on $\{y_3=0\}$ by two relative hyperplanes
   $\{a'y_1+b'y_2=0\}$ and $\{c'y_1+d'y_2=0\}$ and $(a_jy_1+c_jy_2)^2 +
   (b_jy_1+d_jy_2)^2 = (a'y_1+b'y_2)(c'y_1+d'y_2)$.

This is impossible for general choice of $a_j,b_j,c_j,d_j$. In fact, take for
simplicity $b_j=c_j=0$, $a_jd_j\neq 0$. Then the only possible formal
decomposition (up to $\C^*$ is $(a_jy_1)^2 + (d_jy_2)^2 =
   (a_jy_1+id_jy_2)(a_jy_1-id_jy_2)$ (here $i=\sqrt{-1}$). But, since
   ``$a_jy_1$'' is a map from $\hol_B(D_1-\eta_1)$ to $V_2$ and
   ``$d_jy_2$'' is a map from $\hol_B(D_2-\eta_2)$ to $V_2$, these
   factors make sense as relative hyperplanes only when
   $\hol_B(D_1-D_2) \cong \hol_B(\eta_1-\eta_2)$, that is not the case.

It follows that $\C \cap \{y_3=0\}$ is irreducible, then $H \cap \C
\cap \{y_3=0\}$ is $0$-dimensional and therefore there are no
further bicanonical curves through $\Crit(\alpha)$ and $h^1(\T_S) \leq 4$.

\item[($\M_{i,2}$)] The difference with the previous cases is that
$\hol_B(D_1-D_2) \cong \hol_B(\eta_1-\eta_2)$, so, setting as above
$b_6=c_6=0$, $a_6d_6\neq 0$ we can obtain that  $H \cap \C \cap
\{y_3=0\}$ contains a curve by taking $H:=\{a_6y_1\pm
id_6y_2=0\}$. But then $H \cap \Delta$ is
$0$-dimensional of length $3$ so $H \cap \C \cap
\{y_3=0\}$ cannot contain  $\Delta \cap \{y_3=0\}$, that has length $6$.
It follows that there are no further bicanonical curves through
$\Crit(\alpha)$ and $h^1(\T_S) \leq 4$.
\item[($\M_{6,1}$, $\M_{3,1}$, $\M_{4,2}$)] $\C \cap
   \{y_3=0\}$ reduces as union of $\{y_2=0\}$ and $\{a_jy_1+b_jy_2=0\}$,
   that are irreducible for $a_j,b_j$ without common zeroes. The first
   component do not intersect $\Delta$, so to find a bicanonical curve
   we need to take $H$ containing $\{a_jy_1+b_jy_2=0\}$. This is
   possible only when $\hol_B(2\cdot 0 -2 \cdot
   p)$ is the trivial bundle.

Since this is not the case for the three families under consideration,
arguing as above there are no further bicanonical curves through
$\Crit(\alpha)$ and $h^1(\T_S) \leq 4$

\item[($\M_{2,3}$)] Arguing exactly as above we find that the only
   possibility to get a further bicanonical curve through
   $\Crit(\alpha)$ is by  choosing $H:=\{a_1y_1+a_2y_2=0\}$. It
   follows that the set of bicanonical curves through $\Crit(\alpha)$
   is a pencil and therefore $h^1(\T_S)\leq 5$.
\end{itemize}
\end{proof}

\begin{proof}[Proof of theorem \ref{main}]
The first statement comes from remark \ref{comp} and proposition
\ref{dimension}. The second statement is proposition \ref{moduli}. The last
statement was shown in proposition \ref{4famiglie!}.
\end{proof}

\begin{oss}
As mentioned in the introduction the biggest family of minimal surfaces with
$K^2=4$, $p_g=q=1$ constructed by Polizzi is a subfamily of $\M_{2,3}$. We can
be more precise, by looking at the properties of these surfaces (that we will
claim without proof, all follow from the description in \cite{poliziotto}).

It is a family of nodal surfaces obtained as quotient of a product of curves
by an action of $\ZZ_{/2\ZZ} \times \ZZ_{/2\ZZ}$. The group is abelian, so
(arguing as in the proof of \cite{PolThesis}, theorem  6.3)
$\alpha_*\omega^n_S$ in a sum of line bundles for each $n \in \NN$. By
proposition \ref{v1dec} their smooth minimal models give a subfamily of
$\M_{2,3} \cup\M_{4,2} \cup\M_{3,1} \cup\M_{6,1}$.

All Polizzi's surfaces have $4$ nodes.
Since each of our families contains a (smooth minimal) surface with
ample canonical class by proposition \ref{4famiglie!}, and Polizzi's
family is irreducible, then it gives a proper subfamily of one of the
components $\M_{2,3}$, $\M_{4,2}$, $\M_{3,1}$, $\M_{6,1}$. Since it has
dimension $4$, by proposition \ref{dimension} it has codimension $1$
in $\M_{2,3}$.

We can be more precise. The $4$ nodes are contained in two fibres of
the Albanese morphism (two on each fibre), fibres that are
$2$-divisible as Weil divisors on the relative canonical model. It
follows that the singular conics of $\C$ are two double lines. By
the equation of $\C$ in table \ref{Delta}, these are exactly the
surfaces for which $b_1=0$.
\end{oss}

\end{document}